\documentclass[12pt, reqno]{amsart}
\usepackage{latexsym}

\newcommand{\haf}{{\frac{1}{2}}}

\newcommand{\T}{{\mathbb T}}

\newcommand{\C}{{\mathbb C}}

\newcommand{\newsection}[1]
{\subsection{#1}\setcounter{theorem}{0} \setcounter{equation}{0}
\par\noindent}

\newtheorem{theorem}{Theorem}

\newtheorem{lemma}[theorem]{Lemma}

\newtheorem{proposition}[theorem]{Proposition}

\newcommand{\eprop}{\end{proposition}}

\newcommand{\Sum}{\displaystyle\sum}

\begin{document}
\thispagestyle{empty}

\noindent {\large {\bf The hole probability for Gaussian random SU(2)\\ polynomials.\hfill\\
by Scott Zrebiec }}\medskip

\bigskip\bigskip
{\bf Abstract} \par We show that for Gaussian random SU(2)
polynomials of a large degree $N$ the probability that there are
no zeros in the disk of radius $r$ is less than $e^{-c_{1,r}
N^2}$, and is also greater than $e^{-c_{2,r} N^2}$. Enroute to
this result, we also derive a more general result: probability
estimates for the event that the number of complex zeros of a
random polynomial of high degree deviates significantly from its
mean.
\newsection{Introduction and Notation}\par
In this paper we compute the hole probability of one of the key
models for random holomorphic functions: random SU(2) polynomials.
Various properties of the zeros of random SU(2) polynomials have
been studied, in particular the zero point correlation functions
have been computed. This is of particular interest in the physics
literature as the zeros describe a random spin state for the
Majorana representation (modulo phase), \cite{Hannay96}. Further
this choice is intuitively pleasing as the zeros are uniformly
distributed on $\C P^1$ (according to the Fubini-Study metric), or
alternatively the average distribution of zeros is invariant under
the SU(2) action on $\C P^1$.
\par A hole refers to the event where a
particular Gaussian random holomorphic function has no zeros in a
given domain where many are expected. The order of the decay of
the Hole probability has been computed in several cases including
for ``flat" complex Gaussian random holomorphic functions on
$\C^1$, \cite{SodinTsirelson05}, using a method which shall be
used here. This work was subsequently refined to cover other large
deviations in the distribution of the zeros sets,
\cite{Krishnapur}, and generalized to $\C^n$, \cite{Zrebiec06}.
Other results compute the hole probability for a class of complex
Gaussian holomorphic functions on the unit disk,
\cite{PeresVirag04}, and provide a weak general estimate for any
one variable complex Gaussian random holomorphic functions,
\cite{Sodin00}. Additionally significant hole probability results
have been discovered for real Gaussian random polynomials,
(\cite{DemboPoonenShaoZeitouni02}, \cite{LiShao02}).
\par
In this work we will consider the class of random polynomials
whose zeros are distributed on $\C P^1$ according to the
Fubini-Study measure. The random functions which will be studied
here are called Gaussian random SU(2) Polynomials and can be
written as:
$$\psi_{\alpha,N}(z) = \sum_{j=0}^N \alpha_j \sqrt{ N
\choose{j}} z^j \eqno{(1)}$$where $\forall j, \ \alpha_j,$  are
independent identically distributed standard complex Gaussian
random variables (mean 0 and
variance 1). %, with $Var(Re(\alpha_j))=Var(Im(\alpha_j))=\haf $.
For these Gaussian random SU(2) polynomials we will be computing
the hole probability in a manner based on that used by Sodin and
Tsirelson to solve the similar problem for Flat random holomorphic
functions on $\C^1$, \cite{SodinTsirelson05}. First, we shall
estimate the decay rate of the probability for a more general
event:
\begin{theorem}
\label{MainCP1} Let $\psi_{\alpha,N}$ be a degree $N$ Gaussian
random SU(2) polynomial,
$$\psi_{\alpha,N} (z)= \Sum_j \alpha_j
\sqrt{N \choose j} z^j,$$ where $\alpha_j$ are independent
identically distributed complex Gaussian random variables, and let
$\Xi_{\alpha,r,N} $ be the number of zeros of $\psi_{\alpha,N}$
that are in the disk of radius $r$.\par For all $\Delta
> 0,$ and $r>0$ there exists $A_{\Delta,r} $ such
that
$$Prob\left(\left\{\left|\Xi_{\alpha,r,N} - \frac{Nr^2}{1+r^2}
\right| \geq \Delta N  \right\} \right) \leq
%e^{-\left(\frac{c(\Delta^2)}{(1+\frac{\Delta^2}{4})(1+ (1+\Delta)
%r + (1+\Delta)^2 r^2 )^2 }+1\right)\cdot (1-\varepsilon) N^2}
e^{-A_{\Delta,r} N^2}.
$$

%Let $\psi_{\alpha,N}$ be an "Elliptic" Gaussian random
%holomorphic function,
%$$\psi_{\alpha_N} (z)= \Sum_j \alpha_j
%\sqrt{N \choose j} z^j,$$ where $\alpha_j$ are independent
%identically distributed complex Gaussian random variables, and let
%$\Xi_{\alpha,r} $ be the number of zeros of $\psi_\alpha$ that are
%in the disk of radius r.\par For all $\delta
%> 0,$ and $r>0$ there exists $N_{\delta,r} $ and $c_{\delta}>0$ such
%that if $N> N_\delta$ then
%$$Prob\left(\left\{\left|\Xi_{\alpha,r} - \frac{Nr^2}{1+r^2}
%\right| \geq \ \delta {\frac{Nr^2}{1+r^2}} \right\} \right) \leq
%e^{-c_\delta N^2}$$
\end{theorem}
On average, $\psi_{\alpha,N}$ should have $\frac{N r^2}{1+r^2}$
zeros in the disk of radius $r$, as is evident in the above
theorem, or by other more elementary means. Theorem \ref{MainCP1}
gives an upper bound of the rate of decay of the hole probability,
and we will be able to prove a lower bound for the decay rate of
the same order:
\begin{theorem}
\label{Hole probability CP1} Let $Hole_{N,r}=\{(\alpha_j)\in
\C^N:\forall z \in B(0,r), \ \psi_{\alpha,N} (z)\neq 0 \},$ then
there exists $ c_{1,r}, \ and \ c_{2,r}
>0$ such that
$$ e^{-c_{2,r} N^2} \leq Prob (Hole_{N,r})\leq e^{-c_{1,r} N^2}$$
\end{theorem}

The techniques of this paper will generalize to higher dimensions,
as in \cite{Zrebiec06}.\par

Random polynomials of the form studied here are the simplest
examples of a class of natural random holomorphic sections of
large $N$ powers of a positive line bundle on a compact K\"{a}hler
manifold. Most of the results stated in this paper may be restated
in terms of Szeg\"{o} kernels, which exhibit universal behavior in
the large $N$ limit. Hopefully, this paper may provide insight
into proving a similar decay rate for this more general setting.
This has already been done for other properties of random
holomorphic sections, e.g. correlation functions,
\cite{BleherShiffmanZelditchUniv}.\par

\par {\bf Acknowledgement:} I would like
to thank Bernard Shiffman for many useful discussions.

\newsection{SU(2) Invariance}\par

We begin by calling the set of polynomials in one variable whose
degree is less than or equal to $N$, $Poly_N$. $Poly_N$ becomes a
Hilbert space with respect to the SU(2) invariant norm, \cite{BleherShiffmanZelditchUniv}, \cite{SodinTsirelsonRandomComplexZerosI}:\\

$$\displaystyle{\| f\|_{N}
^2}\displaystyle{:=}\displaystyle{\frac{N+1}{\pi}\int_{z\in\C}
|f(z)|^2 \frac{dm(z)}{(1+|z|^2)^{N+2}}},$$ where $dm$ is just the
usual Lebesque measure. For this norm $ \left\{ \sqrt{N \choose j}
z^j\right\}$ is an orthonormal basis, as is $ \left\{ \sqrt{N
\choose j} (az+b)^j (-\overline{b}z+a)^{N-j}\right\}$, where
$|a|^2+|b|^2=1$. Specifically, one orthonormal basis which will be
useful in subsequent work is for any $\zeta\in \C $, $\left\{
\frac{\sqrt{N\choose j} (z - \zeta)^j(1+\overline{\zeta} z)^{N-j}
}{(1+|\zeta|^2)^{\frac{N}{2}}}\right\}$.

Clearly, by line (1), a Gaussian random SU(2) polynomial is
defined as, $\psi_{\alpha,N}(z)= \Sum_{j=0}^{j=N} \alpha_j
\psi_j(z)$, where $\alpha_j$ are i.i.d. standard complex Gaussian
random variables, and $\{\psi_j\}$ is a particular orthonormal
basis. Any basis for $Poly_N$ could have been used and the
Gaussian random SU(2) polynomials would be probabilistically
identical, \cite{SodinTsirelsonRandomComplexZerosI}, as for $
(\alpha_0, \alpha_1, \ldots, \alpha_N )\in \C^{N+1}$ there exists
$U$ a unitary matrix such that for $(\alpha')^T = U \cdot\alpha^T$
$$\psi_{\alpha,N}(z)= \sum \alpha_j \sqrt{N \choose j} z^j =  \frac{\sum \alpha'_j \sqrt{N \choose j} (z-\zeta)^j(1+ \overline{\zeta} z)^{N-j} }{(1+|\zeta|^2)^{\frac{N}{2}}}. \eqno{(2)} $$

%\begin{proof}
%As the set $\displaystyle{\left\{ \frac{\sqrt{\N \choose j}
%(z-\zeta)^j(1+\overline{\zeta}z)^{N-j}}{(|1+\overline{\zeta}z|^2+\left|z-\zeta\right|^2)^\frac{N}{2}}
%\right\}}$ is an orthonormal basis this result is immediate
%Gaussian random holomorphic functions are well defined independent
%of Basis, \cite{Sodin00}.
%\end{proof}

%This result is proven in \cite{SodinTsirelsonRandomComplexZerosI}
%by looking at the SU(2) symmetries. It is also implicitly shown
%there that the probability for an event concerning the zero set of
%random SU(2) polynomials in $B(0,r)$ is identical to the same
%property for $\C\backslash B(0,\frac{1}{r})$.

\newsection{Large deviations of the Maximum of a random SU(2) polynomial}\par

We will use following elementary estimates to compute upper and
lower bounds for the probability of several events:

\begin{proposition}\label{Gauss}
Let $\alpha$ be a standard complex Gaussian random variable, \\
\begin{tabular}{cl}
then & i) $Prob(\{| \alpha| \geq \lambda \}) =
e^{-\lambda^2}$\\
& ii) $Prob(\{ | \alpha| \leq \lambda \}) =1- e^{-\lambda^2}\in
[\frac{\lambda^2}{2}, \lambda^2], if \lambda \leq 1$
\end{tabular}

\end{proposition}

%\begin{proof}
%$\Rightarrow_i$ follows by definition\\
%$\Rightarrow_{ii} \nu(\{ \alpha: |\alpha|\leq
%\lambda\})=1-e^{\lambda^2}= \lambda^2 -\frac{\lambda^4}{2!}+
%\ldots \in [\frac{\lambda^2}{2}, \lambda^2] $
%\end{proof}
This next lemma is key as it states that the maximum of the norm
of a random SU(2) polynomial on the disk of radius $r$ tends to
not be too far from its expected value.
\begin{lemma}\label{GrowthRateCP1}\label{GrowthRateCP1NR0}
For all $\delta\in (0,1]$, and for all $r>0$ there exists
$a_{r,\delta}>0 $ such that\\
%\begin{tabular}{rl}
%1)& $Prob\left( \{ \max_{B(0,r)} |\psi_{\alpha,N}(z)| >
% (1+r^2)^{\frac{N}{2}} (1+\delta)^{\frac{N}{2}}\}\right)<
%e^{-N^2 (1-\varepsilon)} $\\

%2)& $Prob\left( \{ \max_{B(0,r)} |\psi_{\alpha,N}(z)| <
%(1+r^2)^{\frac{N}{2}}(1-\delta)^{\frac{N}{2}}\}\right)< e^{\frac{-
%(2^{\frac{1}{4}} -1) N^2}{4}
%\frac{\delta^2}{1+r^2} (1-\varepsilon)} $\\

 $$Prob\left( \{ \max_{B(0,r)} |\psi_{\alpha,N}(z)|
\notin \left[(1+r^2)^{\frac{N}{2}}(1-\delta)^{\frac{N}{2}} ,
(1+r^2)^{\frac{N}{2}} (1+\delta)^{\frac{N}{2}}\right] \}\right)<
e^{- a_{r,\delta} N^2}$$

%b) For all $\delta >0$ and for all $r\in (0,1]$, there exists
%$N_\delta, \ c>0 $ such that for all $N>N_0$,
%$\displaystyle{Prob \left( \{ \max_{B(0,1)} |\psi_{\alpha,N}(z)|<
%(1+r^2)^{\frac{N}{2}}(1-\delta)^{\frac{N}{2}}\}\right)<e^{-cN^2}}$$
%is a small family of events.\par
\end{lemma}
%While a stronger result than lemma \ref{GrowthRateCP1}-b is true
%this is the one and only result in this paper which is sensitive
%to which radius is chosen. This is caused as the techniques used
%work by transferring the information given in the event to
%information on the range of possible values that a specific random
%coefficients may take on, and for which specific coefficients for
%which this works is highly dependent on r. In practice it is very
%easy to do for any specific r. TO BE CHANGED.

\begin{proof} First we shall prove that for all $\varepsilon>0$, $$Prob\left( \{ \max_{B(0,r)} |\psi_{\alpha,N}(z)|
>(1+r^2)^{\frac{N}{2}} (1+\delta)^{\frac{N}{2}} \}\right)<
e^{- (1-\varepsilon)N^2}.$$ To do this we consider the event
$\Omega_N:= \{ \forall j, \ |\alpha_j| \leq N\},$ the complement
of which has probability $\leq (N+1)e^{-N^2}$, by Proposition
\ref{Gauss}.
\par For $\alpha \in \Omega_N$ we will now estimate
$\displaystyle{\max_{B(0,r)}} |\psi_{\alpha,N}|$:\par
\begin{tabular}{rrl} $\displaystyle{\max_{z\in B(0,r)}|\psi_{\alpha,N}(z )|}$ &$\displaystyle{=}$&$\displaystyle{\max_{z\in B(0,r)} \left|\sum \alpha_j {N \choose j}^\haf (z)^j\right|}$ \\
&$\displaystyle{\leq}$ &$\displaystyle{\max_{z\in B(0,r)} \sum |\alpha_j| {N \choose j}^\haf |z|^j}$ \\
&$\displaystyle{\leq}$&$\displaystyle{ \max_{z\in B(0,r)} N
\sqrt{N+1} (1+ |z|^2)^{\frac{N}{2}}}$, by the Schwartz
inequality.\\& $\displaystyle{ =}$& $\displaystyle{ N
\sqrt{N+1}(1+r^2)^{\frac{N}{2}}}$\\
\end{tabular}\par

For all $\delta>0, \lim_{N\rightarrow \infty} \frac{N \sqrt{N+1}}{
(1+\delta)^{\frac{N}{2}}}=0$, therefore there exists $ N_\delta>0$
such that if $N>N_\delta$ then $\frac{N \sqrt{N+1}}{
(1+\delta)^{\frac{N}{2}}}<1$. \par
Hence if $N>N_\delta$ and if
$\displaystyle{\max_{B(0,r)} |\psi_{\alpha,N}(z)|
>(1+r^2)^{\frac{N}{2}}(1+\delta)^{\frac{N}{2}}}$ then $\alpha\in \Omega_N^c $.
This then guarantees that the event of all such $\alpha$ is a
subset of $\Omega_N^c$ and thus has probability equal to
 $(N+1)e^{-N^2}< e^{-(1-\varepsilon)N^2}$, for large $N$.
 This decay rate is independent of $\delta$ and $r$.
\par We complete the proof by showing that:

$$Prob\left( \{ \max_{B(0,r)} |\psi_{\alpha,N}(z)|
< (1+r^2)^{\frac{N}{2}}(1-\delta)^{\frac{N}{2}} \}\right)< e^{-
a_{r,\delta}N^2 }.$$

Consider the event where $\displaystyle{\max_{B(0,r)}
|\psi_{\alpha,N}(z)|=M < (1+r^2)^{\frac{N}{2}}
(1-\delta)^{\frac{N}{2}} }$. The Cauchy estimates for a
holomorphic function state that: $\displaystyle{
|\psi_{\alpha,N}^{(j)}(0)|\leq j! \frac{M}{r^j}}$. Differentiating
equation (1) yields that $\sqrt{ N \choose{j}} j! \alpha_j =
\psi_{\alpha,N}^{(j)}(0)$. Combining this with Stirling's formula
($\sqrt{2\pi j}j^j e^{-j} <j! < \sqrt{2\pi j}j^j e^{-j}
e^{\frac{1}{12}} $) we get that:
\par
\begin{tabular}{rll}$\displaystyle{ |\alpha_j |}$&$ \displaystyle{
\leq}$&$\displaystyle{ \frac{(1+r^2)^{\frac{N}{2}} (1-\delta
)^\frac{N}{2}}{r^j \sqrt{ N \choose{j}}}}$ \\ & $\displaystyle{
\leq}$ & $\displaystyle{e^{\frac{1}{12}}
\frac{(1+r^2)^{\frac{N}{2}} (1-\delta)^{\frac{N}{2}} \sqrt{2\pi}
(j)^{\frac{j+\haf}{2}} (N-j)^{ \haf (N-j +\haf)} }
{r^j N^{\frac{N + \haf}{2}}}}$\\
&$\displaystyle{ \leq}$&$\displaystyle{ (e^{\frac{1}{12}}
\sqrt{2\pi N} ) \cdot \left(\frac{(1+r^2)^{\frac{N}{2}}
(1-\delta)^{\frac{N}{2}} (j)^{\frac{j}{2}} (N-j)^{ \haf (N-j)} }
{r^j N^{\frac{N}{2}}}\right)}$\\
\end{tabular}\\ For the time being we focus on the term in the second parenthesis which we call $A$.
Writing $j$ as $j= x N, \ x\in (0,1)$, and we now have:
$$A = (1-\delta)^\frac{N}{2}\left( \frac{(1+r^2)}{r^{2x }}  \left(\frac{x}{1-x} \right)^x
(1-x)\right)^\frac{N}{2}
$$

If $x=\frac{r^2}{1+r^2}$ then $A= (1-\delta)^\frac{N}{2}$, which
inspires the next
useful estimate:\\
Claim: Let $m_r=\min\left\{ \frac{(2^\frac{1}{4}-1)r^2}{1+r^2},
\frac{(2^\frac{1}{4}-1)}{1+r^2}\right\}$.\\
If $x\in\left[\frac{r^2}{1+r^2}-m_r \delta, \frac{r^2}{1+r^2}+m_r
\delta \right]$ then $(1+r^2) \left(\frac{x}{r^{2} (1-x)}\right)^x
(1-x)< (1-\delta)^{-\frac{1}{4}}$.\par Note: that
$x\in\left[\frac{r^2}{1+r^2}-m_r \delta, \frac{r^2}{1+r^2}+m_r
\delta \right]\subset (0,1)$\\
Proof: Using the concavity of $x^\frac{1}{4}$, a Taylor's series
and that $\delta <1$, we see that
$\left(\frac{1}{1-\delta}\right)^\frac{1}{4}\geq (2^{\frac{1}{4}}
-1) \delta + 1$.\par We then set $x= (1+\Delta)
\frac{r^2}{1+r^2}$, therefore $$\Delta \in
\left[-\min\left\{\frac{(2^\frac{1}{4}-1)}{r^2},
(2^\frac{1}{4}-1)\right\}\delta,
\min\left\{\frac{(2^\frac{1}{4}-1)}{r^2},
(2^\frac{1}{4}-1)\right\} \delta \right].$$ Thus, $1-x =
\frac{1-\Delta r^2}{ 1+ r^2}$, and from this we compute that:
\par
\begin{tabular}{rcl} $\displaystyle{(1+r^2)
\left(\frac{x}{r^{2}
(1-x)}\right)^{x} (1-x)}$&$\displaystyle{=}$&$\displaystyle{(1-\Delta r^2) \left(\frac{1+\Delta }{1-\Delta r^2}\right)^{x} }$\\
& $\displaystyle{=}$&$\displaystyle{ (1+ \Delta )^x (1-\Delta r^2)^{1-x}}$\\
&$\displaystyle{\leq}$& $\displaystyle{ \max\{1 - \Delta r^2, 1+
\Delta \}},$ as $x\in(0,1)$.
\\ & $\displaystyle{\leq}$& $\displaystyle{
1+(2^{\frac{1}{4}} -1) \delta}$
\end{tabular}\\ Proving the claim. \par
%For $m_r =\min\{ \frac{(2^\frac{1}{4} -1)r^2}{1+r^2},
%\frac{(2^\frac{1}{4} -1)}{1+r^2}\} $,
%$\left[\frac{r^2}{1+r^2}-\delta m_r, \frac{r^2}{1+r^2}+ \delta m_r
%\right]  \subset$
%\\ $\left[\frac{r^2}{1+r^2}-\frac{(2^{\frac{1}{4}} -1)
%\delta}{1+r^2}, \frac{r^2}{1+r^2}+\frac{(2^{\frac{1}{4}} -1)
%\delta r^2}{1+r^2} \right]\bigcap (0,1)$.\par

Therefore for $x\in [\frac{r^2}{1+r^2}-m_r\delta,
\frac{r^2}{1+r^2}+ m_r\delta]$, $A < (1-\delta)^{\frac{3N}{8}}$.
This then in turn guarantees that $|\alpha_j|< e^{\frac{1}{12}}
\sqrt{2\pi N} (1-\delta)^{\frac{3N}{8}}$. The probability this
occurs for a single $\alpha_j$ is less than or equal to
$\left(e^{\frac{1}{12}} \sqrt{2\pi N}
(1-\delta)^{\frac{3N}{8}}\right)^2$. Thus the chance it occurs for
all $\alpha_j$, $j\in \left[(\frac{r^2}{1+r^2}-m_r \delta) N,
(\frac{r^2}{1+r^2}+m_r \delta )N \right]$
%\\ $\bigcap
%\{1, 2, \ldots ,N \}$
, is less than or equal to
$$\left(e^{\frac{1}{12}} \sqrt{2\pi N}
(1-\delta)^{\frac{3N}{8}}\right)^{4 \lfloor N  m_r \delta \rfloor
}.$$ The result follows as there exists $a_{r,\delta}$ such that

$$\left(e^{\frac{1}{12}} \sqrt{2\pi N}
(1-\delta)^{\frac{3N}{8}}\right)^{4 \lfloor N  m_r \delta \rfloor
}
%<e^{\frac{- 3 \lfloor m_r \delta N \rfloor \delta N }{2}  }
<e^{-a_{r,\delta}N^2}$$
\end{proof}

%\begin{remark}
%In the above Lemma for all $\delta<1$ it suffices to use any $$\
%a_{r,\delta}< \min \left\{ \frac{3(2^{\frac{1}{4}}
%-1)\delta^2}{2(1+r^2)}, \frac{3(2^{\frac{1}{4}} -1)r^2
%\delta^2}{2(1+r^2)}\right\}.$$
%\end{remark}
A nice application of this lemma, along with line (2) is the
following:

\begin{lemma}\label{GrowthRateCP1General}\label{ValueEstimateCP1}
For all $\Delta \in (0, 1],$ and $\zeta \in \C\backslash\{0 \} $
there exists $N_{\Delta,|\zeta|}$ and $c_{\Delta,\zeta}>0 $, such
that if $N>N_{\Delta, |\zeta|}$ then
$$\displaystyle{Prob(\{ \max_{z\in B(0,
\Delta )} |\psi_{\alpha,N}(z-\zeta)| < (1+
|\zeta|^2)^{\frac{N}{2}} (1-
\Delta)^{\frac{N}{2}}\})<e^{-c_{\Delta, \zeta} N^2}}$$\par

\end{lemma}
\begin{proof} Let $\Delta
<1$, and set $\delta=\frac{\Delta}{2(1+ |\zeta|+|\zeta|^2)}$.
\\ %(2^{\frac{1}{4}} -1) N^2}{4} \frac{\delta^2}{1+r^2} (1-\varepsilon)

$\displaystyle{(1-\delta)^{\frac{N}{2}} }\ \leq \
\displaystyle{\max_{\partial B(0,\delta)} \frac{\left| \Sum
\alpha'_j \sqrt{N \choose j}(z-\zeta)^j
(1+\overline{\zeta}z)^{N-j}\right|}{(1+
|\zeta|^2)^{\frac{N}{2}}(1+\delta^2)^{\frac{N}{2}}} }$, by Lemma
\ref{GrowthRateCP1NR0} and line (2), except for an event whose
probability is less than $e^{- c_{\Delta,\zeta} N^2}$.\par Let
$\phi(z)=\frac{z-\zeta}{1+\overline{\zeta}z}$, so that we may
rewrite the previous equations as:\par

\begin{tabular}{rrl}
$\displaystyle{(1-\delta)^{\frac{N}{2}} }$&$\leq$&$\displaystyle{\left(\max_{\partial B(0,\delta)} \frac{|1+\overline{\zeta} z|^N}{(1+ |\zeta|^2)^{\frac{N}{2}}(1+\delta^2)^{\frac{N}{2}}} \right) \left(\max_{ B(0,\delta)}\left| \psi_{\alpha',N}(\phi(z))\right|\right) }$\\

&$\leq$&$\displaystyle{\left( \frac{(1+|\zeta| \delta)^N}{((1+
|\zeta|)^2(1+\delta^2))^{\frac{N}{2}}} \right)
\left(\max_{B(-\zeta,(2+2
|\zeta|^2)\delta)} \left| \psi_{\alpha',N}(z)\right|\right)}$, since the\\
\end{tabular}\\
image of $\phi\mid_{B(0,\delta)}\subset B(-\zeta, (2+2
|\zeta|^2)\delta)$.

%, by the following computation:\par
%\begin{tabular}{rrl}
%$\displaystyle{\max_{z\in \partial B(0,\delta)} \left|
%\frac{z-\zeta}{1+ \overline{\zeta}z} + \zeta\right|}$&
%$\displaystyle{=}$&$\displaystyle{ \max_{z\in \partial
%B(0,\delta)} \left| \frac{z-\zeta +\zeta +
%z|\zeta|^2}{1+ \overline{\zeta}z} \right| }$\\
%& $\displaystyle{=}$&$\displaystyle{ \delta \max_{z\in \partial B(0,\delta)} \left| \frac{ (1 +|\zeta|^2)}{1+ \overline{\zeta}z} \right| }$\\
%& $\displaystyle{\leq}$&$\displaystyle{ \delta \  \frac{| \zeta^2
%+1| }{1-\delta|\zeta|}}$\\
%& $\displaystyle{\leq}$&$\displaystyle{ 2 \delta (1+ |\zeta|^2)}$\\
%\end{tabular}\\

Rearranging the previous sets of equations we get the result:\\
%, (as $\{\alpha_j'\}$ are i.i.d. standard Gaussian random Variables):\\
\begin{tabular}{rrl}
$\displaystyle{\max_{ B(-\zeta,\Delta)} |\psi_{\alpha',N}(z)| } $
%&$\displaystyle{\geq }$&$\displaystyle{\max_{ B(-\zeta,(2+2
%|\zeta|^2)\delta)} |\psi_{\alpha',N}(z)|}$ \\

& $\displaystyle{\geq}$&$\displaystyle{
\frac{(1+|\zeta|^2)^{\frac{N}{2}}(1+\delta^2)^{\frac{N}{2}} }{
{(1+|\zeta| \delta)^N}  }
\cdot (1-\delta)^{\frac{N}{2}} } $ \\

%& $\displaystyle{\geq}$&$\displaystyle{
%(1+|\zeta|^2)^{\frac{N}{2}}  (1-|\zeta| \delta)^N
% (1-\delta)^{\frac{N}{2}} } $ \\

%& $\displaystyle{\geq}$&$\displaystyle{
%(1+|\zeta|^2)^{\frac{N}{2}}  (1- 2|\zeta| \delta)^{\frac{N}{2}}
% (1-\delta)^{\frac{N}{2}} } $ \\

& $\displaystyle{\geq}$&$\displaystyle{
(1+|\zeta|^2)^{\frac{N}{2}}  (1- (2+2|\zeta|)
\delta)^{\frac{N}{2}}
 } $ \\

%& $\displaystyle{=}$&$\displaystyle{(1+|\zeta|^2)^{\frac{N}{2}} \left(1- \frac{\Delta (2+ 2|\zeta|)}{2 (1+ |\zeta|+ |\zeta|^2)}\right)^\frac{N}{2} }$\\
&$\displaystyle{\geq}$&$\displaystyle{
(1-\Delta)^{\frac{N}{2}}(1+|\zeta|^2)^{\frac{N}{2}}} $
\end{tabular}\\
\end{proof}

%\begin{remark}
%In particular, for Lemma \ref{GrowthRateCP1General}, if $\Delta<1$
%then it suffices to use

%$$\displaystyle{c_{\Delta,\zeta}= \frac{ 3(2^{\frac{1}{4}} -1)
% \Delta^4}{128 (1+|\zeta|+|\zeta|^2)^2 ((1+|\zeta|+|\zeta|^2)^2+\Delta^2)} }$$

%COMMENT= 3.3 with $\delta=\frac{\Delta}{2(1+|\zeta|+|\zeta|^2)^2}
%, \ r=\delta<1$

%\end{remark}

\newsection{Second key lemma}\par
The goal of this section will be to estimate $\int \log
|\psi_{\alpha,N}(re^{i\theta})|\frac{d\theta}{2\pi}$, which will
be accomplished when we prove lemma
\ref{SurfaceIntegralEstimateOnS1}, using the same techniques as in
\cite{SodinTsirelson05}. As $\log(x)$ becomes unbounded near 0, we
will first prove a deviation result for the event where the $L^1$
norm of $\log|\psi_{\alpha,N}|$ is significantly larger than its
max on the same region.
\begin{lemma} \label{oakleyCP1}For all $r>0$ there exists $N_1$ and $c_r>0$ such that for all $N>N_1$,
$$Prob\left\{\int_{\theta \in \T} |\log(|\psi_{\alpha,N} (r
e^{i\theta})|)| \frac{d \theta}{2\pi} \leq
 5 N  \log\left((2) (1+r^2)\right)\right\}<e^{-c_r N^2}$$
\end{lemma}

\begin{proof}
By Lemma \ref{GrowthRateCP1}, $\exists N_1$ such that if $N>N_1$
then, with the exception of an event, $\Omega$, whose probability
is less than $e^{-c_rN^2}$, there exists $\zeta_0 \in
\partial B(0,\frac{1}{2}r)$ such that $\log(|\psi_{\alpha,N} (\zeta_0)|) >
0$. This also implies that:
$$\int_{\theta = 0}^{2\pi} P_r(\zeta_0, r e^{i\theta}) \log(|\psi_{\alpha,N}(r e^{i\theta})|) \frac{d\theta}{2\pi}\geq \log(|\psi(\zeta_0)|)\geq
0,$$Where $P_r$ is the Poisson kernel: $P_r(\zeta, z)=
\frac{r^2-|\zeta^2|}{|z-\zeta|^2}$. Hence,
$$\int_{\theta = 0}^{2\pi} P_r(\zeta_0, r e^{i\theta}) \log^-(|\psi_{\alpha,N}(r
e^{i\theta})|)\frac{d\theta}{2\pi} \leq \int_{\theta = 0}^{2\pi}
P(\zeta_0, r e^{i\theta}) \log^+(|\psi_{\alpha,N}(r
e^{i\theta})|)\frac{d\theta}{2\pi}$$
\par
Now given the event where
$\displaystyle{\log\max_{B(0,r)}}|\psi_\alpha(z)|< \frac{N}{2}
\log \left((2) (1+r^2)\right)$, (whose complement for $N>N_\delta$
has probability less than $e^{-(1-\varepsilon)N^2}$), we see that
$$\displaystyle{ \int_{\theta= 0}^{2\pi} \log^+(|\psi_{\alpha,N}(re^{i \theta})|) \frac{d \theta}{2
\pi} \leq \frac{N}{2} \log \left((2) (1+r^2)\right)}.$$\par Since
$\zeta_0 \in
\partial B(0, \haf r)$ and $z = r e^{i\theta}$, we have:  $\haf r
\leq |z- \zeta_0 | \leq \frac{3}{2} r $. Hence, by using the
formula for the Poisson Kernel, $\displaystyle{\frac{1}{3} \leq
P(\zeta, z)\leq 3 }$. Putting the pieces together
proves the result:\\
$\displaystyle{\int_{\theta =0}^{2\pi} P_r(\zeta_0, r
e^{i\theta})\log^+}(|\psi_{\alpha,N}(r e^{i\theta})|) \frac{d
\theta}{2
\pi}\displaystyle{\leq \frac{3}{2} N \log\left( 2(1+r^2)\right) }$\\
\begin{tabular}{rrl}
$\displaystyle{\int_{\theta = 0}^{2\pi} \log^-}(|\psi_{\alpha,N}(r e^{i\theta})|)\frac{d \theta}{2\pi}$&$\displaystyle{\leq}$& $\displaystyle{ \frac{1}{\min P(\zeta_0, z)} \int_{\theta = 0}^{2\pi} P_r(\zeta_0, r e^{i\theta})\log^+(|\psi_{\alpha,N}|) \frac{d\theta}{2 \pi} }$\\
&
$\displaystyle{\leq}$& $\displaystyle{ 3  \int_{\theta =0}^{2\pi} P_r(\zeta_0, r e^{i\theta})\log^+(|\psi_{\alpha,N}(r e^{i\theta})|) \frac{d\theta}{2\pi}} $\\
& $\displaystyle{\leq}$& $\displaystyle{
 \frac{9}{2} N \log\left(2 (1+r^2)\right)} $ \\
%& $\displaystyle{\leq 42 N \log\left( (1+\delta) (2)\right)} $
\end{tabular}
 \end{proof}

We now arrive at them main result of this section:
\begin{lemma}\label{SurfaceIntegralEstimateOnS1} For all $\Delta\in (0, 1)$, and
for all $r>0$ there exists $c_{\Delta,r}>0 $ such that
$$Prob\left( \left\{ \int_{\theta =0}^{2\pi}
\log(|\psi_{\alpha,N}(r e^{i\theta})|) \frac{d \theta}{2 \pi} <
\frac{N}{2} \log \left((1+r^2) (1-\Delta) \right) \right\}\right)<
e^{-c_{\Delta,r} N^2}
$$
\end{lemma}
%This lemma provides a probability estimate for how often the
%surface integral deviates below the mean.
\begin{proof}

%This proof is very similar to the one used by Sodin and Tsirelson,
%\cite{SodinTsirelson05}, and as such we merely provide an
%outline.\par

With out loss of generality let $\Delta<\Delta_{0,r}$. Set
$\displaystyle{\delta= \frac{\Delta^4}{81}<1}$. Let $m= \lceil
\frac{1}{\delta} \rceil$, and let
$\displaystyle{\kappa=1-\delta^{\frac{1}{4}}}= 1-\Delta$.\par We
partition the circle of radius $\kappa r$, into m disjoint even
length segments. We choose a $\zeta_j$ within $\delta r$ of the
midpoint of each of these segments such that
$$\log (|\psi_{\alpha,N}(\zeta_j)|) >  \frac{N}{2} \log \left((1 +  \kappa^2 r^2) (1-\delta r)\right), \eqno{(3)}$$
which by Lemma \ref{GrowthRateCP1General} may be done  except for
an event whose probability is less than $ e^{-c_{\delta,r} N^2}$.
Therefore there exists $N_\Delta$ such that if $N>N_\Delta $ the
union of these m events has probability less than or equal to
$(\lceil \frac{2}{\Delta}\rceil+1) e^{-c_{\delta,r} N^2} <
e^{-{c'_\delta,r}N^2}$. We now turn to investigating the average
of $\log |\psi_\alpha(z)|$ on the
unit circle by using Riemann integration and line (3):\\
\\
\begin{tabular}{cl}
$\displaystyle{ \frac{N}{2} \log( (1 + \haf \kappa^2 r^2)}$ &
$\displaystyle{ (1-\delta))\leq \  \frac{1}{m} \Sum_{j=1}^{j=m}
\log|\psi_{\omega,N}(\zeta_j)| }$\\& $ \displaystyle{\leq \
\int_{\theta =0}^{2\pi} \left( \Sum_j \frac{1}{m} P_r(\zeta_j, r
e^{i\theta}) \log(|\psi_{\omega,N} (e^{i \theta})|) \frac{d
\theta}{2 \pi}\right) }$
\\
& $\displaystyle{= \ \int_{\theta =0}^{2\pi} \left(\Sum_j
\frac{1}{m}
(P_r(\zeta_j, re^{i\theta}) -1) \right) \log(|\psi_{\omega,N}( r e^{i \theta})|) \frac{d \theta}{2 \pi}}$\\
&$\displaystyle{ \ \ \ \  + \int_{\theta =0}^{2\pi}
\log(|\psi_{\omega,N} (r e^{i\theta})|) \frac{d \theta}{2 \pi}}$
\\
\end{tabular}\\
This will simplify to:

\begin{tabular}{rl}$\int \log(|\psi_{\alpha,N}(re^{i\theta})|)
\frac{d \theta}{2 \pi}\geq$&$\frac{N}{2}\log \left((1 + \kappa^2
r^2 ) (1- \delta r)\right)$\\ &$ -  ( \int |
\log|\psi_{\alpha,N}(re^{i\theta})|| \frac{d \theta}{2 \pi})
\displaystyle{\max_{\theta\in \T}} |\sum_j \frac{1}{m}
(P_r(\zeta_j, re^{i\theta})-1) |$
\end{tabular}
\\

In \cite{SodinTsirelson05}, it was computed that in exactly this situation that:\\
$$\displaystyle{ \max_{ \theta \in [0, 2\pi]} \left|  \sum_j \frac{1}{m}
(P_r(\zeta_j, re^{i\theta}) -1 )\right| \leq C \delta^{\frac{1}{2
}}}$$

Hence, except for an event of probability $<e^{-cN^2}$, by lemma \ref{oakleyCP1}:\\
\begin{tabular}{rcl}$\int \log(|\psi_{\alpha,N}|) \frac{d
\theta}{2 \pi}$ & $\geq$ & $ \frac{N}{2}\log\left((1 + \kappa^2
r^2) (1- \delta r)\right) - 7 N
\log \left( 2(1+r^2) \right) \cdot C \delta^{\frac{1}{2}},$\\

&$\displaystyle{ \geq}$&$\displaystyle{ \frac{N}{2}\log((1+r^2)
(1-3\delta^\frac{1}{4}))}$.
\end{tabular}\\
\end{proof}

%\begin{remark}
%For the previous lemma, when $\Delta<1$, it suffices to use
%$$c_{\delta,r}<\left(\frac{3(2^\frac{1}{4}-1)(\Delta^{12})}{128 (1+r+r^2)^2((1+r+ r^2)^2+1)}\right)$$
%COMMENT: 3.5 w/ $\Delta= \Delta^3, \ |\zeta|=r$
%\end{remark}

\newsection{Main Results}\par
We now have all the tools needed to prove the two main results of
this paper, starting with theorem \ref{MainCP1}.\\

\begin{proof} (of theorem \ref{MainCP1}). Let $\Xi_{\alpha,r,N}= \int_{B(0,r)} Z_{\psi_{\alpha,N}}=$ the number of
zeros inside the disk of radius $r$ counted with multiplicity. We
ignore the null event where there is a zero on the boundary for a
particular Gaussian random SU(2) polynomial.\par
 It suffices to prove the result for small $\Delta$. Let $\delta= \frac{\Delta^2}{4}<1 $. Let $\kappa = 1 +
 \sqrt{\delta}= 1+ \frac{\Delta}{2}
 $.
Let $a_j$ denote the zeros of a fixed $\psi_{\alpha,N}$ that are
inside $B(0,r)$, and $b_j$ those in $B(0,\kappa r)\backslash
B(0,r)$.\par

We start the proof by recalling Jensen's formula:
$$\log |\psi_{\alpha,N}(0)| = -\Sum \log\left(\frac{\kappa r}{|b_j|}\right) -\Sum \log\left(\frac{\kappa r}{|a_j|}\right)+ \int_{\theta =0}^{2\pi} \log| \psi_{\alpha,N} (\kappa r e^{i\theta})| \frac{d\theta}{2\pi} $$
$$\log | \psi_{\alpha,N} (0)| =  -\Sum \log\left(\frac{r}{|a_j|}\right)+ \int_{\theta =0}^{2\pi} \log| \psi_{\alpha,N} ( r e^{i\theta})| \frac{d\theta}{2\pi} $$
Except for an event of probability $\displaystyle{\leq e^{-a N^2}}$, after subtracting the second line from the first we get that, \\
\begin{tabular}{rl}
$\displaystyle{ \Xi_{\alpha,r,N} \log(\kappa)}$ &$\displaystyle{\leq \int_{\theta =0}^{2\pi} \log| \psi_{\alpha,N} ( \kappa r e^{i\theta})| \frac{d\theta}{2\pi} - \int_{\theta =0}^{2\pi} \log| \psi_{\alpha,N} (r e^{i\theta})| \frac{d\theta}{2\pi}}$\\
& $\displaystyle{ \leq \frac{N}{2} \left( \log \left( (1+ \kappa^2
r^2) (1+\delta)\right) - \int_{\theta =0}^{2\pi} \log |\psi_\alpha
( r e^{i\theta})| \frac{d \theta}{2\pi} \right) }$,\\& by Lemma
\ref{GrowthRateCP1}.
\\ & $\displaystyle{ \leq
\frac{N}{2} \left( \log \left( (1+ \kappa^2 r^2 ) (1+\delta)
\right) - \frac{N}{2} \log \left( (1+r^2)
(1-\delta)\right)\right)}$,\\&   by Lemma
\ref{SurfaceIntegralEstimateOnS1}.\\
%& $\displaystyle{ \leq
%\frac{N}{2}  \log \left(1+ \frac{2\sqrt{\delta} r^2+ \delta r^2 +
%2 \delta+ 2\delta r^2 + O (\delta^2)}{(1+r^2)
%} \right)  }$, \\
& $\displaystyle{ \leq \frac{N}{2}   \left( \frac{2\sqrt{\delta}
r^2+ \delta r^2 + 2 \delta +2\delta r^2  }{(1+r^2) } - \frac{2
\delta
r^4}{(1+r^2)^2} + O (\delta^\frac{3}{2}) \right) }$, \\
\end{tabular}\par
Therefore,\\
\begin{tabular}{rcl}$\displaystyle{\Xi_{\alpha,r,N}}$&$\leq$&$\displaystyle{  N  \left( \frac{
r^2+ \haf \sqrt\delta r^2 + \sqrt\delta + \sqrt{\delta}r^2
}{(1+r^2) } - \frac{ \sqrt \delta r^4}{(1+r^2)^2} + O (\delta)
\right)\left( 1 + \frac{\sqrt\delta}{2} +
O(\delta)\right) }$\\
&$\leq$&$\displaystyle{   \frac{N r^2}{1+r^2}+ 2 N \sqrt
\delta + O (\delta) }$\\

\end{tabular}

%Putting the pieces together,\\ \begin{tabular}{rrl}
%$\Xi_{\alpha,r,N}$&$\displaystyle{ \leq}$&$\displaystyle{
%\frac{\frac{N}{2} \left( \frac{2\sqrt{\delta} r^2}{1+r^2} + 2
%\delta - \frac{2r^4 \delta}{(1+r^2)^2} + \frac{\delta r^2}{1+r^2}
%+O(\delta^\frac{3}{2})   \right)} {\sqrt \delta - \frac{\delta}{2}+

%O(\delta^{\frac{3}{2}})}}$\\
%&$\displaystyle{ =}$&$\displaystyle{  \left( \frac{N r^2}{1+r^2} +
%N \sqrt{\delta} - \frac{Nr^4 \sqrt{\delta}}{(1+r^2)^2} + \frac{N
%\sqrt{\delta} r^2}{2(1+r^2)} + O(\delta)   \right) \left(1+
%\frac{\sqrt{\delta}}{2} + O(\delta)\right)}$\\

%&$\displaystyle{=}$&$\displaystyle{\frac{N r^2}{1+r^2} +
% 3\sqrt{\delta} N, \ \forall \delta<\delta_0}$
%\end{tabular}\par

This proves the probability estimate when the number of zeros in
the lower hemisphere is significantly larger than expected. If we
choose new random variables by taking the unitary transformation
which reverses the order of the random variables, $\alpha=
(\alpha_0, \alpha_1,\ldots, \alpha_N), U \alpha = (\alpha_N,
\alpha_{N-1}, \ldots, \alpha_0)= \alpha'$, and assemble a random
SU(2) polynomial from these Gaussian Random variables,
$\psi_{\alpha',N(z)}$, a direct computation shows that for
$z_0\neq0, \  \psi_{\alpha,N}(z_0)= 0 \newline \Leftrightarrow
\psi_{\alpha',N}(\frac{1}{z_0})= 0 $. Thus,
$$\displaystyle{\left\{\alpha':\Xi_{\alpha',
\frac{1}{r}}-\frac{N}{1+r^2}< \Delta N\right\} =
\left\{\alpha:\Xi_{\alpha,r}-\frac{Nr^2}{1+r^2}\geq \Delta
N\right\}}.$$ As, $\{ \alpha_j'\}$ are all i.i.d. Gaussian random
events we may apply the first part of this proof implies that
these two events have probability less than or equal to
$e^{-aN^2}$, and the result follows immediately.

\end{proof}

%\begin{remark}
%In the preceding theorem for $r=1$ it suffices to use
%$$A_{\Delta,r}<\left(\frac{(2^\frac{1}{4}-1)(\Delta^{12})}{3,840}\right)$$
%COMMENT: min of 4.5 and 4.5 with $\frac{1}{r}$
%\end{remark}

We have just implicitly proven an upper bound on the order of the
decay of the hole probability. We will now compute the lower bound
to finish the proof theorem \ref{Hole probability CP1}

\begin{proof} (of theorem \ref{Hole probability CP1}). If $\forall z\in B(0,r), \ \psi_{\alpha,N}(z)\neq 0$ then
$\Xi_{\alpha,r}=0$
, therefore by the previous theorem, $Hole_{N,r}$ is contained in
the event:
$$\left\{\left|\Xi_{\alpha,r} - \frac{Nr^2}{1+r^2}
\right| \geq  \frac{r^2}{1+r^2} \frac{N}{2}  \right\}.
$$
In other words events of the type whose probability has just been
given an upper bound in theorem 1.1.
% specifically setting $\Delta=\frac{1}{7(1+r^2)}$ in the previous theorem.
For this event there exists an $N_r$, such that for all $N>N_r$,
the probability is less than $e^{-a N^2}$.
\par We must still prove the lower bound, for which we start by
considering the event, $\Omega$ which consists of $\alpha_j$
where:\par\begin{tabular}{cl} & $\displaystyle{|\alpha_0|\geq N}$\\
&$\displaystyle{|\alpha_j|< {N \choose j}^\frac{-1}{2} r^{-j}}$
\end{tabular}\\
If $\alpha\in \Omega$, then  $|\alpha_0| > \sum_{j>0} |\alpha_j|
{N \choose j}^\frac{1}{2} r^j$. Hence $\forall z \in B(0,r)\
\psi_{\alpha,N}(z)\neq 0
 \Rightarrow \Omega \subset Hole_{N,r} $.

We will now approximate the probability of $\Omega$ in order to
get a lower bound for the probability of $Hole_{N,r}$:\\
\begin{tabular}{rrl}$\displaystyle{Prob\left(\left\{ \ |\alpha_j|<{N \choose
j}^\frac{-1}{2} r^{-j} \right\}
\right)}$&$\displaystyle{\geq}$&$\displaystyle{
\frac{1}{2} \frac{1}{{N \choose j} r^{2j}}}$, if ${N\choose j}\geq\frac{1}{r^{2j}} $\\
&$\displaystyle{\geq}$&$\displaystyle{ \frac{1}{2} \frac{(N-j)!
j!}{N!} r^{-2j}}$\\
&$\displaystyle{\geq}$&$\displaystyle{\frac{\sqrt{\pi N}
}{2^{N+\frac{3}{2}}} r^{-2j}e^{-\frac{1}{12}}}$
\\ & $\displaystyle{\geq}$&$\displaystyle{
\frac{e^{\frac{-1}{12}}}{(2)^{N}} \min\left\{ \frac{1}{r},
\frac{1}{r^{2N}} \right\}, \ if \ N>1 \ and \ j\geq
1}$\\&$\displaystyle{=}$& $\displaystyle{ e^{-N\log(2) -
\max\{\log(r), 2N\log(r) \} -\frac{1}{12} }}$
\end{tabular}\\
$Prob(\{|\alpha_0| > N \})= e^{-N^2}$\\
\begin{tabular}{rrl}Hence, $\displaystyle{Prob (\Omega)}$&$\displaystyle{\geq}$&$\displaystyle{ e^{-N^2-  N^2 \log(2) - N \max\{\log(r), 2N\log(r) \}- \frac{N}{12} }}$\\
&$\displaystyle{\geq}$&$ \displaystyle{ e^{-c_r N^2}}$
\end{tabular}
\end{proof}

%\begin{remark}
%The hole probability for the unit disk obeys the following
%inequalities:
%$$ e^{-(\frac{11}{12}+ \frac{\log(2)}{2})N^2}< Prob(Hole_{N,1})< e^{-\left(\frac{N^2}{186,026,899,981,440}\right)}$$
%COMMENT:$c_{2,r} = \frac{11}{12}+ \frac{\log(2)}{2}+ \log(r)$, 5.1
%w/ $\Delta=\frac{1}{7(1+r^2)} $.
%\end{remark}

\newsection{Generalizations to m dimensions.}
\par In this paper we have just computed the decay rate of the
hole probability for Gaussian random SU(2) polynomials. This work
may be adapted to the case of m variables or random SU(m+1)
polynomials: $$\psi_{\alpha,N,m}(z_1, \ldots, z_m) = \Sum_{j=
(j_1,\ldots, j_m)} \alpha_j \sqrt{N \choose j_1, \ldots, \ j_m }
z^j
$$where $\alpha_j$ are i.i.d. standard complex Gaussian random
variables. For this system, by following the steps in this paper
and in \cite{Zrebiec06}, the following theorem will be proven in
subsequent work.
\begin{theorem}
Let $Hole_{N,r,m}$ be the event where for all $z\in B(0,r)$,
\newline $\psi_{\alpha, N,m}\neq 0$. Then there exists $N_{r,m},\ \ c_{r,m}>0$ such that for all $N>N_{r,m}$
$$e^{-c_{r,m}N^{m+1}}\leq Prob(Hole_{N,r,m})\leq
e^{-c'_{r,m}N^{m+1}}$$
\end{theorem}
Also, a result similar to theorem \ref{MainCP1}, but with an
``$e^{-cN^{m+1}}$" rate of decay could be proven.

\end{document}